\documentclass[11pt]{article}
\usepackage{amssymb}
\newtheorem{lemma}{Lemma}[section]
\newtheorem{definition}[lemma]{Definition}

\newtheorem{theorem}[lemma]{Theorem}
\newtheorem{proposition}[lemma]{Proposition}
\newtheorem{corollary}[lemma]{Corollary}
\newtheorem{example}[lemma]{Example}

\newtheorem{note}[lemma]{Remark}

\def\endproof{\hfill$\Box$}

\def\endproof{\hfill$\Box$}

\title{On freedom and independence \\ in hypergraphs of models of theories\footnote{This research was partially
supported by Committee of Science in Education and Science
Ministry of the Republic of Kazakhstan (Grant No. AP05132546) and
Russian Foundation for Basic Researches (Project No.
17-01-00531-a).}}

\author{}
\author{B.Sh.~Kulpeshov, S.V.~Sudoplatov}
\date{}

\begin{document}

\maketitle

\begin{abstract}
Notions of freedom and independence for hypergraphs of models of a
theory are defined. Properties of these notions and their
applications to some natural classes of theories are studied.
\end{abstract}


\bigskip
Hypergraphs of models of a theory are related to derived objects,
allowing to obtain an essential structural information on both
theories themselves and related semantical objects including graph
ones \cite{CCMCT14, Su013, Su08, Baik, SudKar16, KulSud17, Sud17,
KulSud171, KS12018}.

In the present paper notions of freedom and independence for
hypergraphs of models of a theory are defined. Properties of these
notions and their applications to some natural classes of theories
are studied.

\section{Preliminaries}
\noindent

Recall that a {\it hypergraph\/} is any pair of sets $(X,Y)$,
where $Y$ is some subset of the Boolean $\mathcal{P}(X)$ of a set
$X$. The set $X$ is called the {\em universe} of the hypergraph
$(X,Y)$, and elements of $Y$ are {\em edges} of the hypergraph
$(X,Y)$.

Let $\mathcal{M}$ be some model of a complete theory $T$.
Following \cite{SudKar16} we denote by $H(\mathcal{M})$ the family
of all subsets $N$ of the universe $M$ of the structure
$\mathcal{M}$, which are universes of elementary submodels
$\mathcal{N}$ of the model $\mathcal{M}$: $H(\mathcal{M})=\{N\mid
\mathcal{N}\preccurlyeq\mathcal{M}\}$. The pair
$(M,H(\mathcal{M}))$ is called the {\it hypergraph of elementary
submodels\/} of the model $\mathcal{M}$ and it is denoted by
$\mathcal{H}(\mathcal{M})$\index{$\mathcal{H}(\mathcal{M})$}.

For a cardinality $\lambda$ we denote by
$H_\lambda(\mathcal{M})$\index{$H_\lambda(\mathcal{M})$} and
$\mathcal{H}_\lambda(\mathcal{M})$\index{$\mathcal{H}_\lambda(\mathcal{M})$},
respectively, the restrictions of $H(\mathcal{M})$ and
$\mathcal{H}(\mathcal{M})$ on the class of elementary submodels
$\mathcal{N}$ of $\mathcal{M}$ such that $|N|<$~$\lambda$.

By
$\mathcal{H}_p(\mathcal{M})$,\index{$\mathcal{H}_p(\mathcal{M})$}
$\mathcal{H}_l(\mathcal{M})$,\index{$\mathcal{H}_l(\mathcal{M})$}
$\mathcal{H}_{\rm npl}(\mathcal{M})$,\index{$\mathcal{H}_{\rm
npl}(\mathcal{M})$}
$\mathcal{H}_h(\mathcal{M})$,\index{$\mathcal{H}_h(\mathcal{M})$}
$\mathcal{H}_s(\mathcal{M})$\index{$\mathcal{H}_s(\mathcal{M})$}
we denote the restrictions of the hyper\-graph
$\mathcal{H}_{\omega_1}(\mathcal{M})$ on the class of elementary
submodels $\mathcal{N}$ of the model $\mathcal{M}$, that are prime
over finite sets, limit, non-prime and non-limit, homogeneous,
saturated respectively. Similarly, by
$H_p(\mathcal{M})$,\index{$H_p(\mathcal{M})$}
$H_l(\mathcal{M})$,\index{$H_l(\mathcal{M})$} $H_{\rm
npl}(\mathcal{M})$,\index{$H_{\rm npl}(\mathcal{M})$}
$H_h(\mathcal{M})$,\index{$H_h(\mathcal{M})$}
$H_s(\mathcal{M})$\index{$H_s(\mathcal{M})$} are denoted the
correspondent restrictions for  $H_{\omega_1}(\mathcal{M})$.

\begin{definition} {\rm \cite{SudKar16, Engel}.
Let $(X,Y)$ be a hypergraph, $x_1,x_2$ be distinct elements of
$X$. We say that the element $x_1$ is {\em separated} or {\em
separable} from the element $x_2$, or {\em $T_0$-separable} if
there is $y\in Y$ such that $x_1\in y$ and $x_2\notin y$. The
elements $x_1$ and $x_2$ are called {\em separable}, {\em
$T_2$-separable}, or {\em Hausdorff separable} if there are
disjoint $y_1,y_2\in Y$ such that $x_1\in y_1$ and $x_2\in y_2$.}
\end{definition}

\begin{theorem}\label{th21} {\rm \cite{SudKar16}.} Let $\mathcal{M}$ be an $\omega$-saturated model of
a countable complete theory $T$, $a$ and $b$ be elements of
$\mathcal{M}$. The following are equivalent:

$(1)$ the element $a$ is separable from the element $b$ in
$\mathcal{H}(\mathcal{M})$;

$(2)$ the element $a$ is separable from the element $b$ in
$\mathcal{H}_{\omega_1}(\mathcal{M})$;

$(3)$ $b\notin{\rm acl}(a)$.
\end{theorem}

\begin{theorem}\label{th26} {\rm \cite{SudKar16}.} Let $\mathcal{M}$ be an $\omega$-saturated model of
a countable complete theory $T$, $a$ and $b$ be elements of
$\mathcal{M}$. The following are equivalent:

$(1)$ the elements $a$ and $b$ are separable in
$\mathcal{H}(\mathcal{M})$;

$(2)$ the elements $a$ and $b$ are separable in
$\mathcal{H}_{\omega_1}(\mathcal{M})$;

$(3)$ ${\rm acl}(a)\cap{\rm acl}(b)=\varnothing$.
\end{theorem}

\begin{corollary}\label{co27} {\rm \cite{SudKar16}.}
Let $\mathcal{M}$ be an $\omega$-saturated model of a countable
complete theory $T$, $a$ and $b$ be elements of $\mathcal{M}$, and
there exists the prime model over $a$. The following are
equivalent:

$(1)$ the element $a$ is separable from the element $b$ in
$\mathcal{H}(\mathcal{M})$;

$(2)$ the element $a$ is separable from the element $b$ in
$\mathcal{H}_{\omega_1}(\mathcal{M})$;

$(3)$ the element $a$ is separable from the element $b$ in
$\mathcal{H}_p(\mathcal{M})$;

$(4)$ $b\notin{\rm acl}(a)$.
\end{corollary}

\begin{corollary}\label{co28} {\rm \cite{SudKar16}.}
Let $\mathcal{M}$ be an $\omega$-saturated model of a countable
complete theory $T$, $a$ and $b$ be elements of $\mathcal{M}$, and
there exist the prime models over $a$ and $b$ respectively. The
following are equivalent:

$(1)$ the elements $a$ and $b$ are separable in
$\mathcal{H}(\mathcal{M})$;

$(2)$ the elements $a$ and $b$ are separable in
$\mathcal{H}_{\omega_1}(\mathcal{M})$;

$(3)$ the elements $a$ and $b$ are separable in
$\mathcal{H}_p(\mathcal{M})$;

$(4)$ ${\rm acl}(a)\cap{\rm acl}(b)=\varnothing$.
\end{corollary}

\begin{definition} {\rm \cite{SudKar16}.
Let $(X,Y)$ be a hypergraph, $X_1,X_2$ be disjoint nonempty
subsets of the set $X$. We say that the set $X_1$ is {\em
separated} or {\em separable} from the set $X_2$, or {\em
$T_0$-separable} if there is $y\in Y$ such that $X_1\subseteq y$
and $X_2\cap y=\varnothing$. The sets $X_1$ and $X_2$ are called
{\em separable}, {\em $T_2$-separable}, or {\em Hausdorff
separable} if there are disjunct $y_1,y_2\in Y$ such that
$X_1\subseteq y_1$ and $X_2\subseteq y_2$.}
\end{definition}

By using proofs of theorems \ref{th21} and \ref{th26}, the
following generalizations of these theorems are established.

\begin{theorem}\label{th31} {\rm \cite{SudKar16}} Let $\mathcal{M}$
be a $\lambda$-saturated model of a complete theory $T$,
$\lambda\geq{\rm max}\{|\Sigma(T)|$, $\omega\}$, $A$ and $B$ be
nonempty sets in $\mathcal{M}$ having the cardinalities
$<\lambda$. The following are equivalent:

$(1)$ the set $A$ is separable from the set $B$ in
$\mathcal{H}(\mathcal{M})$;

$(2)$ the set $A$ is separable from the set $B$ in
$\mathcal{H}_\lambda(\mathcal{M})$;

$(3)$ ${\rm acl}(A)\cap B=\varnothing$.
\end{theorem}

\begin{theorem}\label{th32} {\rm \cite{SudKar16}} Let $\mathcal{M}$ be a $\lambda$-saturated model of a complete
theory $T$,  $\lambda\geq{\rm max}\{|\Sigma(T)|$, $\omega\}$, $A$
è $B$ be nonempty sets in $\mathcal{M}$ having the cardinalities
$<\lambda$. The following are equivalent:

$(1)$ the sets $A$ and $B$ are separable in
$\mathcal{H}(\mathcal{M})$;

$(2)$ the sets $A$ and $B$ are separable in
$\mathcal{H}_\lambda(\mathcal{M})$;

$(3)$ ${\rm acl}(A)\cap{\rm acl}(B)=\varnothing$.
\end{theorem}

We obtain by analogy with corollaries \ref{co27} and \ref{co28}

\begin{corollary}\label{co33} {\rm \cite{SudKar16}.}
Let $\mathcal{M}$ be an $\omega$-saturated model of a small theory
$T$, $A$ and $B$ be finite nonempty sets in $\mathcal{M}$. The
following are equivalent:

$(1)$ the set $A$ is separable from the set $B$ in
$\mathcal{H}(\mathcal{M})$;

$(2)$ the set $A$ is separable from the set $B$ in
$\mathcal{H}_{\omega_1}(\mathcal{M})$;

$(3)$ the set $A$ is separable from the set $B$ in
$\mathcal{H}_p(\mathcal{M})$;

$(4)$ ${\rm acl}(A)\cap B=\varnothing$.
\end{corollary}

\begin{corollary}\label{co34} {\rm \cite{SudKar16}.}
Let $\mathcal{M}$ be an $\omega$-saturated model of a small theory
$T$, $A$ and $B$ be finite nonempty sets in $\mathcal{M}$. The
following are equivalent:

$(1)$ the sets $A$ and $B$ are separable in
$\mathcal{H}(\mathcal{M})$;

$(2)$ the sets $A$ and $B$ are separable in
$\mathcal{H}_{\omega_1}(\mathcal{M})$;

$(3)$ the sets $A$ and $B$ are separable in
$\mathcal{H}_p(\mathcal{M})$;

$(4)$ ${\rm acl}(A)\cap{\rm acl}(B)=\varnothing$.
\end{corollary}

\begin{definition} {\rm \cite{KulSud171}.
Let $\mathcal{M}$ be a model of a theory $T$ with a hypergraph
$\mathcal{H}=(M,H(\mathcal{M}))$ of elementary submodels, $A$ be
an infinite definable set in $\mathcal{M}$, of arity $n$:
$A\subseteq M^n$. The set $A$ is called {\em $\mathcal{H}$-free}
if for any infinite set $A'\subseteq A$, $A'=A\cap Z^n$ for some
$Z\in H(\mathcal{M})$ containing parameters for $A$. Two
$\mathcal{H}$-free sets $A$ and $B$ of arities $m$ and $n$
respectively are called {\em $\mathcal{H}$-independent} if for any
infinite $A'\subseteq A$ and $B'\subseteq B$ there is $Z\in
H(\mathcal{M})$ containing parameters for $A$ and $B$ and such
that $A'=A\cap Z^m$ and $B'=B\cap Z^n$.}
\end{definition}

Note the following properties \cite{KulSud171}.

1. Any two tuples of a $\mathcal{H}$-free set $A$, whose distinct
tuples do not have common coordinates, have same type.

Indeed, if there are tuples $\bar{a},\bar{b}\in A$ with ${\rm
tp}(\bar{a})\ne{\rm tp}(\bar{b})$ then for some formula
$\varphi(\bar{x})$ the sets of solutions of that formula and of
the formula $\neg\varphi(\bar{x})$ divide the set $A$ into two
nonempty parts $A_1$ and $A_2$, where at least one part, say
$A_1$, is infinite. Taking $A_1$ for $A'$ we have $A_1=A\cap Z^n$
for appropriate $Z\in H(\mathcal{M})$ and~$n$. Then by the
condition for tuples in $A$ we have $A_2\cap Z^n=\emptyset$ that
is impossible since $Z$ is the universe of an elementary submodel
of $\mathcal{M}$.

Thus the formula $\varphi(\bar{x})$, defining $A$, implies some
complete type in $S^n(\emptyset)$, and if $A$ is
$\emptyset$-definable then $\varphi(\bar{x})$ is a principal
formula.

In particular, if the set $A$ is $\mathcal{H}$-free and
$A\subseteq M$, then the formula, defining $A$, implies some
complete type in $S^1(\emptyset)$.

2. If $A\subseteq M$ is a $\mathcal{H}$-free set, then $A$ does
not have nontrivial definable subsets, with parameters in $A$,
i.e., subsets distinct to subsets defined by equalities and
inequalities with elements in $A$.

Indeed, if $B\subset A$ is a nontrivial definable subset then $B$
is defined by a tuple $\bar{a}$ of parameters in $A$, forming a
{\em finite} set $A_0\subset A$, and $B$ is distinct to subsets of
$A_0$ and to $A\setminus C$, where $C\subseteq A_0$. Then removing
from $A$ a set $B\setminus A_0$ or $(A\setminus B)\setminus A_0$,
we obtain some $Z\in H(\mathcal{M})$ violating the satisfiability
for $B$ or its complement. It contradicts the condition that $Z$
is the universe of an elementary submode of $\mathcal{M}$.

3. If $A$ and $B$ are two $\mathcal{H}$-independent sets, where
$A\cup B$ does not have distinct tuples with common coordinates,
then $A\cap B=\emptyset$.

Indeed, if $A\cap B$ contains a tuple $\bar{a}$, then, choosing
infinite sets $A'\subseteq A$ and $B'\subseteq B$ with $\bar{a}\in
A'$ and $\bar{a}\notin B'$, we obtain $\bar{a}\in A'=A\cap Z^n$
for appropriate $Z\in H(\mathcal{M})$ and $n$, as so $\bar{a}\in
B\cap Z^n=B'$. This contradiction means that $A\cap B=\emptyset$.

\begin{definition} {\rm \cite{KulSud17}. The {\em complete union}
of hypergraphs $(X_i,Y_i)$, $i\in I$, is the hypergraph
$\left(\bigcup\limits_{i\in I}X_i,Y\right)$, where
$Y=\left\{\bigcup\limits_{i\in I}Z_i\mid Z_i\in Y_i\right\}$. If
the sets $X_i$ are disjoint, the complete union is called {\em
disjoint} too.  If the set $X_i$ form a $\subseteq$-chain, then
the complete union is called {\em chain}.}
\end{definition}

By Property 3 we have the following theorem on decomposition of
restrictions of hyper\-graphs $\mathcal{H}$, representable by unions
of families of $\mathcal{H}$-independent sets.

\begin{theorem}\label{dcu} {\rm \cite{KulSud171}.}
A restriction of hypergraph $\mathcal{H}=(M,H(\mathcal{M}))$ to a
union of a family of $\mathcal{H}$-free $\mathcal{H}$-independent
sets $A_i\subseteq M$ is represented as a disjoint complete union
of restrictions $\mathcal{H}_i$ of the hypergraph $\mathcal{H}$ to
the sets $A_i$.
\end{theorem}

Proof. Consider a family of $\mathcal{H}$-independent sets
$A_i\subseteq M$. By Property 3 these sets are disjoint, and using
the definition of $\mathcal{H}$-independence we immediately obtain
that the union of restrictions $\mathcal{H}_i$ of $\mathcal{H}$ to
the sets $A_i$ is complete.
\endproof

\medskip
\begin{definition} {\rm \cite{KS2}. Let $\mathcal{M}$ be some model of a complete theory $T$,
$(M,H(\mathcal{M}))$ be a hypergraph of elementary submodels of the model $\mathcal{M}$.
Sets $N\in H(\mathcal{M})$ are called {\em elementarily submodel} or {\em elementarily
substructural} in $\mathcal{M}$.}
\end{definition}

\begin{proposition}\label{press1} {\rm \cite{KS2}}
Let $A$ be a definable set in an $\omega_1$-saturated model $\mathcal{M}$ of
a countable complete theory $T$. Then exactly one of the following conditions holds:

$(1)$ The set $A$ is finite and is contained in any elementarily
substructural set in $\mathcal{M}$;

$(2)$ The set $A$ is infinite, has infinitely many different
intersections with elementarily substructural sets in
$\mathcal{M}$, and all these intersections are infinite; and the
indicated intersections can be chosen so as to form an infinite
chain/antichain by inclusion.
\end{proposition}

\begin{proposition}\label{press2} {\rm \cite{KS2}}
Let $A$ be a definable set in the countable saturated model $\mathcal{M}$ of a small theory $T$.
Then exactly of the following conditions holds:

$(1)$ The set $A$ is finite and is contained in any elementarily
substructural set in $\mathcal{M}$;

$(2)$ The set $A$ is infinite, has infinitely many different
intersections with elementarily substructural sets in
$\mathcal{M}$, and all these intersections are infinite; and the
indicated intersections can be chosen so as to form an infinite
chain / antichain by inclusion.
\end{proposition}

\medskip
Note that the above concepts and statements by a natural manner are transferred to
hypergraphs  $\mathcal{H}_\lambda(\mathcal{M})$, $\mathcal{H}_p(\mathcal{M})$,
$\mathcal{H}_l(\mathcal{M})$, $\mathcal{H}_{\rm
npl}(\mathcal{M})$, $\mathcal{H}_h(\mathcal{M})$,
$\mathcal{H}_s(\mathcal{M})$.

\vskip 3mm

Recall that a subset $A$ of a linearly ordered structure
$\mathcal{M}$ is called {\it convex} if for any $a, b\in A$ and
$c\in M$ whenever $a<c<b$ we have $c\in A$. A {\it weakly
o-minimal structure} is a linearly ordered structure
$\mathcal{M}=\langle M,=,<,\ldots \rangle$ such that any definable
(with parameters) subset of the structure $\mathcal{M}$ is a union
of finitely many convex sets in $\mathcal{M}$.

In the following definitions $\mathcal{M}$ is a weakly o-minimal
structure, $A, B\subseteq M$, $\mathcal{M}$ be $|A|^+$-saturated,
$p,q\in S_1(A)$ be non-algebraic types.
\begin{definition}\rm \cite{bbs1}
We say that $p$ is not {\it weakly orthogonal} to $q$
($p\not\perp^w  q$) if there exist an $A$-definable formula
$H(x,y)$, $\alpha \in p(\mathcal{M})$ and $\beta_1, \beta_2 \in
q(\mathcal{M})$ such that $\beta_1 \in H(\mathcal{M},\alpha)$ and
$\beta_2 \not\in H(\mathcal{M},\alpha)$.
\end{definition}

\begin{definition} \rm \cite{k2003} We say that $p$ is not {\it quite orthogonal} to
 $q$ ($p\not\perp^q q$) if there exists an $A$-definable bijection $f: p(\mathcal{M})\to q(\mathcal{M})$.
 We say that a weakly o-minimal theory is {\it quite o-minimal} if the notions of
 weak and quite orthogonality of 1-types coincide.
\end{definition}

In the paper \cite{KS} the countable spectrum for quite o-minimal
theories with non-maximal number of countable models has been
described:

\begin{theorem}\label{KS_apal} Let $T$ be a quite o-minimal theory with non-maximum many
countable models. Then $T$ has exactly $3^k\cdot 6^s$ countable
models, where $k$ and $s$ are natural numbers. Moreover, for any
$k,s\in\omega$ there exists a quite o-minimal theory $T$ having
exactly $3^k\cdot 6^s$ countable models.
\end{theorem}

Realizations of these theories with a finitely many countable
models are natural ge\-ne\-ra\-li\-za\-tions of Ehrenfeucht
examples obtained by expansions of dense linear orderings by a
countable set of constants, and they are called theories of  {\em
Ehrenfeucht type}. Moreover, these realizations are representative
examples for hypergraphs of prime models \cite{CCMCT14, Su08,
SudKar16}.  We consider operators for hypergraphs allowing on one
hand to describe the decomposition of hypergraphs of prime models
for quite o-minimal theories with few countable models, and on the
other hand pointing out constructions leading to the building of
required hypergraphs by some simplest ones.

Denote by $(M,H_{\rm dlo}(\mathcal{M}))$ a hypergraph of (prime)
elementary submodels of a countable model $\mathcal{M}$ of the theory
of dense linear order without endpoints.

\begin{note}\label{no441} {\rm
The class of hypergraphs $(M,H_{\rm dlo}(\mathcal{M}))$ is closed
under countable chain complete unions, modulo density and having
an encompassing dense linear order without endpoints. Thus, any
hypergraph $(M,H_{\rm dlo}(\mathcal{M}))$ is represented as a
countable chain complete, modulo density, union of some its proper
subhypergraphs.  }
\end{note}

Any countable model of a theory of Ehrenfeucht type is a disjoint
union of some intervals, which are ordered both themselves and
between them, and of some singletons. Dense subsets of the
intervals form universes of elementary substructures. So, in view
of Remark \ref{no441}, we have:

\begin{theorem}\label{KS_apr} {\rm \cite{KulSud17}.}
A hypergraph of prime models of a countable model of a theory of
Ehrenfeucht type is represented as a disjoint complete, modulo
density, union of some hypergraphs in the form $(M,H_{\rm
dlo}(\mathcal{M}))$ as well as singleton hypergraphs of the form
$(\{c\},\{\{c\}\})$.
\end{theorem}

\begin{note}\label{no442} {\rm Taking into consideration links
between sets of realizations of $1$-types, which are not weakly
orthogonal, as well as definable equivalence relations, the
construction for the proof of Theorem \ref{KS_apr} admits a
natural generalization for an arbitrary quite o-minimal theory
with few countable models. Here conditional complete unions should
be additionally {\em coordinated}, i.e., considering definable
bijections between sets of realizations of $1$-types, which are
not quite orthogonal.}
\end{note}

\section{On relative freedom and independence in hypergraphs of models of theories}
\noindent

As shown in Section 1, $\mathcal{H}$-free sets does not have
non-trivial definable subsets. By this note at studying subsets
$A'$ of definable sets $A\subseteq M^n$ in structures
$\mathcal{M}$ of a non-empty signature, where $A'=A\cap (M_1)^n$
for some $\mathcal{M}_1\prec\mathcal{M}$, it is naturally instead
of ``absolute'' $\mathcal{H}$-freedom to consider relative
$\mathcal{H}$-freedom taking into account, as for dense linear
orders, the specifics of subsets $A'$ by some syntactical
information taken from the complete diagram $D^\ast(\mathcal{M})$
of the system $\mathcal{M}$. In the following section we take into
account this specific for ordered theories, and in this section we
introduce general notions of relative $\mathcal{H}$-freedom and
$\mathcal{H}$-independence, and we also establish links between
distinct types of relativity.

\begin{definition}
{\rm Let $\mathcal{M}$ be some model of a theory $T$ with a
hypergraph of elementary submodels
$\mathcal{H}=(M,H(\mathcal{M}))$, $D^\ast(\mathcal{M})$ be the
complete diagram of the model $\mathcal{M}$, ${\bf D}$ be some set
of diagrams $\Phi(A_0)\subseteq D^\ast(\mathcal{M})$ such that for
some language $\Sigma\subseteq\Sigma(\mathcal{M})$ if
$\varphi(\bar{a})$ is a quantifier-free formula of the language
$\Sigma$, $\bar{a}\in A_0$, then $\varphi(\bar{a})\in\Phi(A_0)$ or
$\neg\varphi(\bar{a})\in\Phi(A_0)$. Here, the set $A_0$ is called
the {\em universe} of the diagram $\Phi(A_0)$. We say that a set
$A\subseteq M^n$ {\em satisfies} a diagram $\Psi\in{\bf D}$ if
$\Psi=\Phi(A_0)$ for the set $A_0$ consisting of all the
coordinates of tuples from $A$. The set $A\subseteq M^n$ is called
{\em relatively $\mathcal{H}$-free}, $\mathcal{H}$-free modulo
${\bf D}$, or {\em $(\mathcal{H},{\bf D})$-free} if for any set
$A'\subseteq A$ satisfying some diagram of ${\bf D}$ the equality
$A'=A\cap Z^n$ holds for some $Z\in H(\mathcal{M})$ containing
parameters for $A$. Two $(\mathcal{H},{\bf D})$-free sets $A$ and
$B$ of arity $m$ and $n$ respectively are called {\em relatively
$\mathcal{H}$-independent}, $\mathcal{H}$-independent modulo ${\bf
D}$, or {\em $(\mathcal{H},{\bf D})$-independent} if for any sets
$A'\subseteq A$ and $B'\subseteq B$ satisfying some diagrams of
${\bf D}$ there exists $Z\in H(\mathcal{M})$ containing parameters
for $A$ and $B$ and such that $A'=A\cap Z^m$ and $B'=B\cap Z^n$.}
\end{definition}

Note that at defining ``absolute'' $\mathcal{H}$-freedom and
$\mathcal{H}$-independence as $\Sigma$ it is considered the empty
language, a set $A$ is definable and infinite, and diagrams $\Phi$
are taken either quantifier-free on all infinite sets $A'\subseteq
A$ or as result of adding to these quantifier-free diagrams
schemes of infinity of sets for $A'$.

Unlike definability in the case of type definability or
non-definability of a set $A$ under consideration of relative
freedom and independence both a scheme of infinity and an infinity
itself of sets $A'$ may not be required.  Indeed, for a theory $T$
of unary predicates $P_i$ with $P_{i+1}\subset P_i$, $i\in\omega$,
the non-isolated type $p(x)=\{P_i(x)\mid i\in\omega\}$ can have
the set of realizations of any, finite or infinite, cardinality.
Thus, the set $A=p(\mathcal{M})$, $\mathcal{M}\models T$,
non-having non-trivial connections is $(\mathcal{H},{\bf
D}_p)$-free for a set of diagrams ${\bf D}_p$ describing
realizability of the type $p(x)$ by elements of an arbitrary set
$A'\subseteq A$.

If the theory $T$ is expanded by unary predicates $Q_i$ with
conditions $Q_{i+1}\subset Q_i\subset\overline{P_0}$,
$i\in\omega$, then the set $B=q(\mathcal{M})$, where
$q(x)=\{Q_i(x)\mid i\in\omega\}$ is free relatively a set of
diagrams ${\bf D}_q$ describing realizability of the type $q(x)$
by elements of an arbitrary set $B'\subseteq B$, will be
$(\mathcal{H},{\bf D}_p\cup{\bf D}_q)$-independent.

\medskip
Since at extending diagrams of ${\bf D}$ a family of considered
sets can only decrease the following hold:

\medskip
{\bf Monotonicity properties.} 1. If ${\bf D}\subseteq{\bf
D}'\subseteq\mathcal{P}(D^\ast(\mathcal{M}))$ and a set
$A$ is $(\mathcal{H},{\bf D}')$-free then
$A$ is $(\mathcal{H},{\bf D})$-free.

2. If ${\bf D}\subseteq{\bf
D}'\subseteq\mathcal{P}(D^\ast(\mathcal{M}))$, sets $A$ and $B$
are $(\mathcal{H},{\bf D}')$-independent then $A$ and $B$ are
$(\mathcal{H},{\bf D})$-independent.

3. If diagrams of ${\bf D}$ are some restrictions/extensions of
suitable diagrams from the set ${\bf
D}'\subseteq\mathcal{P}(D^\ast(\mathcal{M}))$, with preservation
of their universes, and the set $A$ is $(\mathcal{H},{\bf
D})$-free, then $A$ is $(\mathcal{H},{\bf D}')$-free.

4. If diagrams of ${\bf D}$ are some restrictions/extensions of
suitable diagrams of the set ${\bf
D}'\subseteq\mathcal{P}(D^\ast(\mathcal{M}))$, with preservation
of their universes, the sets $A$ and $B$ are $(\mathcal{H},{\bf
D})$-independent, then $A$ and $B$ are $(\mathcal{H},{\bf
D}')$-independent.

\medskip
Inverse implications in monotonicity properties are not true.
Indeed, if an infinite definable set $A\subseteq M$ is partitioned
by a unary predicate $P$ into two non-empty parts then $A$ is not
$\mathcal{H}$-free, although this set can be considered as
$(\mathcal{H},{\bf D}')$-free, for the language $\{P^{(1)}\}$,
where ${\bf D}'$ consists of diagrams describing cardinalities
$|A\cap P|$ and $|A\cap\overline{P}|$. The considered effect, at
which two disjoint infinite definable sets $A$ and $B$ are
partitioned by the predicate $P$ into non-empty disjoint parts,
shows that an independence of the sets $A$ and $B$ can be failed
at transition from $(\mathcal{H},{\bf D}')$ to $(\mathcal{H},{\bf
D})$.

Further for simplicity we will mostly consider the notions of
relative freedom and in\-de\-pen\-dence for sets $A\subseteq M$,
although these considerations can be adapted, for example, by the
operation $\mathcal{M}^{\rm eq}$, for arbitrary sets $A\subseteq
M^n$.

In connection with the introduced concepts, a series of natural
questions and problems arises.

\medskip
1. Is an arbitrary set in the given structure free relative to
some set of diagrams ${\bf D}$?

\medskip
2. Characterize the condition of $(\mathcal{H},{\bf D})$-freedom of a set.

\medskip
3. Characterize the condition of $(\mathcal{H},{\bf D})$-independence of sets.

\medskip
4. Is there and if yes then which, a condition on sets of diagrams
${\bf D}'$ such that sets $A$ are $(\mathcal{H},{\bf D}')$-free,
but not $(\mathcal{H},{\bf D})$-free when ${\bf D}\subset{\bf
D}'$?

\medskip
5. Is there and if yes then which, a condition on sets of diagrams
${\bf D}'$ such that sets $A$ and $B$ are $(\mathcal{H},{\bf
D}')$-independent, but not $(\mathcal{H},{\bf D})$-independent
when ${\bf D}\subset{\bf D}'$?

\medskip
One of the ways of answering these questions is the considered
below choice for sets $A$ of diagrams $\Phi(A_0)$ with suitable
sets $A_0$. However, this approach does not take into account the
structural specificity of sets $A$ and thus not in the full extent
it reflects the real freedom of these sets, as well as their
independence. This specificity is taken into account in the
following sections, for some specific classes of theories.

\begin{proposition}\label{f1}
For any set $A\subseteq M$ in a model $\mathcal{M}$ of a theory
$T$ there exists the set of diagrams ${\bf D}$ such that $A$ is
$(\mathcal{H},{\bf D})$-free.
\end{proposition}

Proof. It is sufficiently to take as ${\bf D}$ an arbitrary set of
diagrams $\Phi\subseteq D^\ast(\mathcal{M})$ of which the
universes contain the set $A$. \endproof

\begin{proposition}\label{f2}
For any sets $A,B\subseteq M$ in a model $\mathcal{M}$ of a theory
$T$ there exists the set of diagrams ${\bf D}$ such that the sets
$A$ and $B$ are $(\mathcal{H},{\bf D})$-independent.
\end{proposition}

Proof.  It is sufficiently to take as ${\bf D}$ an arbitrary set
of diagrams $\Phi\subseteq D^\ast(\mathcal{M})$ of which the
universes contain the set $A\cup B$. \endproof

\medskip
By Propositions \ref{press1} and \ref{press2} the following statement holds.

\begin{proposition}\label{f3}
For any set $A\subseteq M$ in an $\omega_1$-saturated {\rm
(}countable saturated{\rm )} model $\mathcal{M}$ of a countable {\rm
(}small{\rm )} theory $T$ exactly one of the following conditions holds:

$(1)$ $A$ is finite and $(\mathcal{H},{\bf D})$-free only relative
to the set of diagrams ${\bf D}$ whose universes contain the set
$A$;

$(2)$ $A$ is infinite and $(\mathcal{H},{\bf D})$-free relative to
some infinite set of diagrams ${\bf D}$ whose universes have
infinite distinct intersections with $A$, and these intersections
can be chosen so that they form an infinite chain/antichain by
inclusion.
\end{proposition}

\section{On freedom and independence in hypergraphs of models of theories
with unary
predicates and theories with equivalence relations}
\noindent

In this section we describe decompositions of hypergraphs
$\mathcal{H}(\mathcal{M})$ for theories with unary predicates and
some theories with equivalence relations.

Firstly we consider the theory $T$ with unary predicates $P_i$,
$i\in I$, and an equivalence relation coinciding with the equality
relation. Since all language connections between ele\-ments are
bounded by the condition of their uniformity, i.e., coincidence of
$1$-types, a description of a hypergraph
$\mathcal{H}(\mathcal{M})$ is reduced to description of its
restrictions on the sets of realizations of complete $1$-types.

Due to the lack of connections between $1$-types, based on
Proposition \ref{press1} the following assertion holds for
definable sets consisting of realizations of isolated types.

\begin{proposition}\label{prun}
For any definable set $A$ consisting of realizations in a model
$\mathcal{M}$ of some isolated $1$-type $p$, either
$|H(\mathcal{M})|=1$ when $|A|<\omega$ or
$|H(\mathcal{M})|=2^\lambda$ when $|A|=\lambda\geq\omega$. Here,
infinite sets $A$ are $\mathcal{H}$-free and
$\mathcal{H}$-independent.
\end{proposition}

By Theorem \ref{dcu} and Proposition \ref{prun} the following holds

\begin{corollary}\label{coun}
The restriction of a hypergraph $\mathcal{H}=(M,H(\mathcal{M}))$
on a union of any family of sets $A_j\subseteq M$, each of which
is the set of realizations of some isolated $1$-type $p_j$, is
represented in the form of disjoint complete union
$\mathcal{H}_{\rm isol}$ of restrictions $\mathcal{H}_j$ of the
hypergraph $\mathcal{H}$ on sets $A_j$.
\end{corollary}

Now we consider restrictions of a hypergraph
$\mathcal{H}(\mathcal{M})$ on the sets $B_k$ of realizations of
non-isolated types $q_k$. Since these types can be both omitted
and realized by an arbitrary quantity of realizations,
restrictions $\mathcal{H}(\mathcal{M})\upharpoonright B_k$ are
represented in the form of atomic Boolean lattices $L_k$.

If considered types $q_k$ can be omitted in aggregate (for
example, if the theory has a prime model) then the family of
lattices $L_k$ compose their complete union also forming an atomic
Boolean lattice:

\begin{proposition}\label{prun2}
Restrictions of a hypergraph $\mathcal{H}=(M,H(\mathcal{M}))$ on a
union of any family of sets $B_k\subseteq M$, each of which is the
set of realizations of some non-isolated $1$-type $q_k$, where the
types $q_k$ can be omitted in aggregate, are represented in the
form of disjoint complete union
$\mathcal{H}_{n\mbox{\scriptsize\rm -isol}}$ of restrictions
$\mathcal{H}_k$ of the hypergraph $\mathcal{H}$ on sets $B_k$.
This disjoint complete union forms an atomic Boolean lattice.
\end{proposition}

At consideration of restriction of a hypergraph
$\mathcal{H}=(M,H(\mathcal{M}))$ on a union of a family of sets
$A_j$ and $B_k$, if types $q_k$ can be omitted, a representation
of this restriction in the form of disjoint complete union of
restrictions $\mathcal{H}_j$ and $\mathcal{H}_k$, and also in the
form of disjoint complete union $\mathcal{H}_{\rm isol}$ and
$\mathcal{H}_{n\mbox{\scriptsize\rm -isol}}$ holds:

\begin{proposition}\label{prun3}
Restriction of a hypergraph $\mathcal{H}=(M,H(\mathcal{M}))$ on a
union of any family of sets $A_j\subseteq M$, each of which os the
set of realizations of some isolated $1$-type $p_j$, and also any
family of sets $B_k\subseteq M$, each of which is the set of
realizations of some non-isolated $1$-type $q_k$, where the types
$q_k$ can be omitted in aggregate, is represented in the form of
disjoint complete union $\mathcal{H}'$ of restrictions
$\mathcal{H}_j$ of the hypergraph $\mathcal{H}$ on sets $A_j$ and
restrictions $\mathcal{H}_k$ of the hypergraph $\mathcal{H}$ on
sets $B_k$, i.e. in the form of disjoint complete union
$\mathcal{H}_{\rm isol}$ and $\mathcal{H}_{n\mbox{\scriptsize\rm
-isol}}$. This disjoint complete union forms an atomic Boolean
lattice modulo $\mathcal{H}_{\rm isol}$.
\end{proposition}

To complete a description of decomposition of a hypergraph
$\mathcal{H}(\mathcal{M})$ it remains to consider its restrictions
on the sets $C_l$ of realizations of non-isolated types $r_l$,
non-omitted in aggregate. Such families of types arise, for
example, in the theory of independent unary predicates, non-having
isolated $1$-types. In this case subsets $C_l$ also can be varied
arbitrarily, but with condition $C$ of satisfaction of all
consistent formulas by elements from $A_j$, $B_k$, $C_l$. Thus, a
{\em conditional} complete union or {\em $C$-union} of
restrictions $\mathcal{H}_l$ of the hypergraph $\mathcal{H}$ on
sets $C_l$ arises:

\begin{proposition}\label{prun4}
Restriction of a hypergraph $\mathcal{H}=(M,H(\mathcal{M}))$ on a
union of any family of sets $A_j\subseteq M$, each of which is the
set of realizations of some isolated $1$-type $p_j$, any family of
sets $B_k\subseteq M$, each of which is the set of realizations of
some non-isolated $1$-type $q_k$, where the types $q_k$ can be
omitted in aggregate, and also any family of sets $C_l\subseteq
M$, each of which is th set of realizations of some non-isolated
$1$-type $r_l$, where the types $r_l$ cannot be omitted in
aggregate, is represented in the form of  $C$-union of
restrictions $\mathcal{H}_j$ of the hypergraph $\mathcal{H}$ on
sets $A_j$, restrictions $\mathcal{H}_k$ of the hypergraph
$\mathcal{H}$ on sets $B_k$, and also restrictions $\mathcal{H}_l$
of the hypergraph $\mathcal{H}$ on sets $C_l$.
\end{proposition}

On the base of statements \ref{prun}, \ref{coun}, \ref{prun2},
\ref{prun3}, \ref{prun4} the following theorem holds describing
the de\-com\-po\-sition of a hypergraph $\mathcal{H}(\mathcal{M})$
by means of four types of hypergraphs holds.

\begin{theorem}\label{thun}
For any model $\mathcal{M}$ of a theory $T$ of some unary
predicates the hypergraph
$\mathcal{H}(\mathcal{M})=(M,H(\mathcal{M}))$ is represented in
the form of disjoint complete union of some of the following
hypergraphs:

$1)$ a hypergraph with the universe $M_0$ consisting of
realizations of all algebraic $1$-types, and having the only edge
coinciding with $M_0$;

$2)$ a disjoint complete union of $\mathcal{H}$-free hypergraphs
of which the universes consist of realizations of non-algebraic
isolated $1$-types;

$3)$ a disjoint complete union of hypergraphs forming atomic
Boolean lattices, the universes of which consist of realizations
of non-isolated $1$-types omitted in aggregate;

$4)$ a $C$-union of hypergraphs the universes of which consist of
realizations of non-isolated $1$-types non-omitted in aggregate.
\end{theorem}

\begin{example} \rm As examples including all the types of hypergraphs 1)--4),
described in Theo\-rem \ref{thun}, it can be considered
hypergraphs of disjoint unions of the following structures:

i) structures consisting of unique non-empty finite unary
predicates;

ii) structures consisting of unique non-empty infinite unary
predicates;

iii) structures consisting of countably many disjoint non-empty
unary predicates;

iv) structures consisting of countably many independent unary
predicates.
\end{example}

Now we consider theories $T$ with equivalence relations $E_i$,
$i\in I$.

If the relation $E_i$ is unique then in the theory there is an
information on the number and cardinality of equivalence classes.
If the number of these classes is finite then all of them are
presented in any elementary submodel $\mathcal{N}$ of a model
$\mathcal{M}\models T$ and a hypergraph
$\mathcal{H}=(M,H(\mathcal{M}))$ is represented in the form of
disjoint complete union of its restrictions on $E_i$-classes. The
same is related to finite $E_i$-classes with a finite number for
given cardinality $n$. If the number of such $E_i$-classes is
infinite then in elementary submodels $\mathcal{N}$ of the model
$\mathcal{M}$ are included arbitrary families of $n$-element
$E_i$-classes, i.e., the hypergraph $\mathcal{H}$ is {\em free}
relative to $E_i$-classes.

At consideration of infinite $E_i$-classes each of which is
$\mathcal{H}$-free and distinct $E_i$-classes are
$\mathcal{H}$-independent. If the number of infinite $E_i$-classes
is finite then each of which is presented in models $\mathcal{N}$,
and if the number of infinite $E_i$-classes is infinite then a
{\em $\mathcal{H}$-freedom} for $E_i$-classes holds, i.e. any
infinite subset of these $E_i$-classes together with given finite
$E_i$-classes form an elementary submodel with the universe of
$H(\mathcal{M})$.

The indicated $\mathcal{H}$-freedom and $\mathcal{H}$-independence
is extended to theories with successively embedded equivalence
relations under the condition of uniformity $E_i$-classes. If in
uniform equivalence classes there exist structures with unary
predicates then at representation of a hypergraph  $\mathcal{H}$
in the form of disjoint complete union of restrictions on these
classes also provides a decomposition described in Theorem
\ref{thun}.

In general case the problem of describing hypergraphs for theories
with equivalence re\-la\-tions remains open.

\section{On freedom and independence in hypergraphs of models of ordered theories}
\noindent

Note that if we consider an arbitrary non-algebraic isolated type
$p\in S_1(\emptyset)$ in an arbitrary almost $\omega$-categorical
quite o-minimal theory $T$ then in any model $\mathcal{M}\models
T$ the set $p(\mathcal{M})$ will not be $\mathcal{H}$-free, since
if we take as $A'$ some closed interval $[a, b]\subset
p(\mathcal{M})$, where $a<b$, then there is no $\mathcal{M}_1\prec
\mathcal{M}$ such that $A'=p(\mathcal{M})\cap M_1$.

Another reason for the violation of $\mathcal{H}$-freedom is the
possibility of taking an infinite set $A'\subset p(\mathcal{M})$
that is not dense, while the sets $p(\mathcal{M})\cap M_1$ for
models of the theory $T$ must be dense.

An arbitrary open interval containing an element $b$ is said to be
a {\it neighbourhood} of the element $b$. Recall that an arbitrary
subset $A$ of a linearly ordered structure $\mathcal{M}$ is {\it
open} if for any $b\in A$ there is a neighbourhood of the element
$b$ that is contained in $A$.

\begin{definition}
{\rm Let $T$ be a weakly o-minimal theory, $\mathcal{M}\models T$,
$p\in S_1(\emptyset)$ be a non-algebraic isolated type. We say
that $p(\mathcal{M})$ is {\em relatively $\mathcal{H}$-free}, {\em
$\mathcal{H}$-free relative to convex sets}, or {\em
$(\mathcal{H}, {\rm cs})$-free} if for any open convex set
$A'\subseteq p(\mathcal{M})$ the equality $A'=p(\mathcal{M})\cap
M_1$ holds for some $M_1\in H(\mathcal{M})$.}
\end{definition}

We note that, in addition to the fact that $\mathcal{H}$
hypergraphs allow to select all infinite subsets of
$\mathcal{H}$-free sets, the corresponding hypergraphs for
$(\mathcal{H}, {\rm cs}) $-free sets in addition to convex sets
without endpoints make it possible to isolate dense sets without
endpoints.  For example, for a theory of dense linear order
without endpoints any dense subset without endpoints is
distinguished by such way.

We also introduce the following necessary definitions.

\begin{definition} \rm  \cite{CCMCT14, IPT}
Let $p_1(x_1),\ldots,p_n(x_n)\in S_1(T)$. A type
$q(x_1,\ldots,x_n)\in S(T)$ is called
\emph{$(p_1,\ldots,p_n)$-type} if
$q(x_1,\ldots,x_n)\supseteq\bigcup\limits_{i=1}^n p_i(x_i)$. The
set of all \ $(p_1,\ldots,p_n)$-types of the theory  $T$ is
denoted by $S_{p_1,\ldots,p_n}(T)$. A countable theory $T$ is
called {\em almost $\omega$-categorical} if for any types
$p_1(x_1),\ldots,p_n(x_n)\in S(T)$ there is only finitely many
types $q(x_1,\ldots,x_n)\in S_{p_1,\ldots,p_n}(T)$.
\end{definition}

\begin{definition}\rm \cite{bbs2}
Let $\mathcal{M}$ be a weakly o-minimal structure, $A\subseteq M$,
$\mathcal{M}$ be $|A|^+$-saturated, $p\in S_1(A)$ be a
non-algebraic type.

(1) An $A$-definable formula $F(x, y)$ is {\it $p$-preserving} or
{\it $p$-stable} if there are $\alpha$, $\gamma_1$, $\gamma_2 \in
p(\mathcal{M})$ such that $F(\mathcal{M},\alpha)\setminus
\{\alpha\}\ne\emptyset$ è $\gamma_1 < F(\mathcal{M}, \alpha) <
\gamma_2$.

(2) A $p$-preserving formula $F(x, y)$ is {\it convex to right
(left)} if there is $\alpha$ $\in$ $p(\mathcal{M})$ such that
$F(\mathcal{M},\alpha)$ is convex, $\alpha$ is the left (right)
endpoint of the set $F(\mathcal{M},\alpha)$ and $\alpha\in
F(\mathcal{M},\alpha)$.
 \end{definition}

\begin{definition}\rm \cite{BaizKu}
We say that a $p$-preserving convex to right (left) formula
$F(x,y)$ is {\it equivalence-generating} if for any $\alpha, \beta
\in p(\mathcal{M})$ such that $\mathcal{M}\models F(\beta,
\alpha)$, the following holds:
$$\mathcal{M}\models \forall x [x\ge \beta \to [F(x,\alpha) \leftrightarrow
F(x,\beta)]]\quad (\mathcal{M}\models \forall x [x\le \beta \to
[F(x,\alpha) \leftrightarrow F(x,\beta)]])$$
\end{definition}

\begin{definition}  \rm \cite{k1}
Let $T$ be a weakly o-minimal theory, $\mathcal{M}$ be a
sufficiently saturated model of $T,$ and let $\phi(x)$ be an
arbitrary $M$-definable formula with one free variable. The {\it
convexity rank of the formula} $\phi(x)$ $(RC(\phi(x)))$ is
defined as follows:

1)  $RC(\phi(x)) \geq 1$ if $\phi(\mathcal{M})$ is infinite.

2) $RC(\phi(x)) \geq  \alpha  +  1$ if there are a parametrically
definable equivalence relation $E(x,y)$ and infinitely many
elements $b_i, i\in \omega $, such that:
\begin{itemize}
\item For any  $i, j\in \omega$ whenever $i\ne j$ we have
$\mathcal{M} \models \neg E(b_i, b_j)$;
\item For each $i\in \omega \quad RC(E(x,  b_i))
\geq \alpha$ and $E(\mathcal{M},b_i)$ is a convex subset of
$\phi(\mathcal{M})$.
\end{itemize}

3) $RC(\phi(x)) \geq \delta$ if $RC(\phi(x)) \geq \alpha$ for all
$\alpha < \delta \; (\delta$ is limit).

If $RC(\phi(x)) = \alpha$ for some $\alpha$ then we say that
$RC(\phi(x))$ is defined. Otherwise (i.e. if $RC(\phi(x))$ $\geq$
$\alpha$ for all $\alpha$), we put $RC(\phi(x)) =\infty$.
\end{definition}

The {\it convexity rank of an $1$-type} $p$ ($RC(p)$) is called
the infimum of the set $\{RC(\phi(x)) \mid \phi(x)\in p\}$, i.e.
$RC(p):= \inf \{RC(\phi(x)) \mid \phi(x)\in p\}$.

\begin{lemma}\label{l_free}
Let $T$ be an almost $\omega$--categorical quite o-minimal theory,
$\mathcal{M}\models T$, $p\in S_1(\emptyset)$ be a non-algebraic
isolated type. Then $p(\mathcal{M})$ is relatively
$\mathcal{H}$-free $\Leftrightarrow$ $RC(p)=1$.
\end{lemma}

Proof of Lemma \ref{l_free}. $(\Rightarrow)$ Let
$\;p(\mathcal{M})$ is relatively $\mathcal{H}$-free. Assume the
contrary: $RC(p)>1$. By binarity of $T$ there is an
$\emptyset$-definable equivalence relation $E(x,y)$ partitioning
$p(\mathcal{M})$ into infinitely many infinite convex sets. Take
an arbitrary $a\in p(\mathcal{M})$ and consider
$E(a,\mathcal{M})$. Obviously, $E(a,\mathcal{M})$ is open convex
set, and there is no an elementary submodel $\mathcal{M}_1$ of
$\mathcal{M}$ such that $E(a,\mathcal{M})=p(\mathcal{M})\cap M_1$.

$(\Leftarrow)$ Let $RC(p)=1$. We argue to show that
$p(\mathcal{M})$ is indiscernible over $\emptyset$. By binarity of
$T$ it is sufficiently to prove that $p(\mathcal{M})$ is
2-indiscernible over $\emptyset$. Assume the contrary: there are
$\langle a_1, a_2\rangle$, $\langle a'_1, a'_2\rangle\in
[p(\mathcal{M})]^2$ such that  $a_1<a_2$, $a'_1<a'_2$ and
$tp(\langle a_1, a_2\rangle/\emptyset)\ne tp(\langle a'_1,
a'_2\rangle /\emptyset)$. Then there exists $a''_2\in
p(\mathcal{M})$ such that $a_1<a''_2$ and $tp(\langle
a_1,a_2\rangle/\emptyset)\ne tp(\langle a_1, a''_2\rangle
/\emptyset)$. Consequently, there is an $\emptyset$-definable
formula $\phi(x,y)$ such that $\mathcal{M}\models
\phi(a_1,a_2)\land \neg \phi(a_1, a''_2)$. By weak o-minimality we
can assume that $\phi(a_1, \mathcal{M})$ is convex. Without loss
of generality, we will also assume that $a_2<a''_2$. Then consider
the following formula:
$$F(x,a_1):=x\ge a_1\land \exists y [\phi(a_1, y)\land x\le y].$$

It is easy to see that $F(x,y)$ is a $p$-preserving convex to
right. If $F(x,y)$ is equivalence-generating, we have a
contradiction with $RC(p)=1$. If $F(x,y)$ is not
equivalence-generating, it contradicts to almost
$\omega$-categoricity of $T$. Thus, $p(\mathcal{M})$ is
indiscernible over $\emptyset$, whence for any open convex set
$A'\subseteq p(\mathcal{M})$ there is an elementary submodel
$\mathcal{M}_1$ of $\mathcal{M}$ such that $A'=p(\mathcal{M})\cap
M_1$.
\endproof

\begin{example}\rm
Let $\mathcal{M}=\langle\mathbb{Q},<,f^1\rangle$ be a linearly
ordered structure, where $\mathbb{Q}$ is the set of rational
numbers, $f(x)=x+1$ is an unary function on $\mathbb{Q}$.

It is easily seen that $\mathcal{M}$ is an o-minimal structure,
and $Th(\mathcal{M})$ is not almost $\omega$-categorical. Note
also that $p(x):=\{x=x\}\in S_1(\emptyset)$ is a non-algebraic
isolated type, $RC(p)=1$, but $p(\mathcal{M})$ is not relatively
$\mathcal{H}$-free.
\end{example}

\begin{definition}
{\rm Let $T$ be a weakly o-minimal theory, $\mathcal{M}\models T$,
$p,q\in S_1(\emptyset)$ be non-algebraic isolated types,
$RC(p)=RC(q)=1$. We say that $p(\mathcal{M})$ and $q(\mathcal{M})$
are {\it relatively $\mathcal{H}$-independent}, {\em
$\mathcal{H}$-independent with regard to convex sets}, or {\em
$(\mathcal{H}, {\rm cs})$-independent} if for any open convex sets
$A'\subseteq p(\mathcal{M})$ and $B'\subseteq q(\mathcal{M})$
there is $M_1\in H(\mathcal{M})$ such that $A'=p(\mathcal{M})\cap
M_1$ and $B'=q(\mathcal{M})\cap M_1$.}
\end{definition}

\vskip 2mm Let  $p_1,p_2,\ldots, p_s\in S_1(\emptyset)$ be
non-algebraic types. We say that the family of types $\{p_1,
\ldots, p_s\}$ is {\it orthogonal over $\emptyset$} if for any
sequence $(n_1$, $\ldots$,$n_s)\in \omega^s$, for any increasing
tuples $\bar a_1, \bar a'_1 \in [p_1(\mathcal{M})]^{n_1}$,
$\ldots$, $\bar a_s, \bar a'_s\in [p_s(\mathcal{M})]^{n_s}$ such
that $tp(\bar a_1/\emptyset)=tp(\bar a'_1/\emptyset)$, $\ldots$,
$tp(\bar a_s/\emptyset)=tp(\bar a'_s/\emptyset)$ we have
$tp(\langle \bar a_1, \ldots$, $\bar
a_s\rangle/\emptyset)=tp(\langle \bar a'_1$, $\ldots$, $\bar
a'_s\rangle /\emptyset)$.

\begin{lemma}\label{l_ind}
Let $T$ be an almost $\omega$-categorical quite o-minimal theory,
$\mathcal{M}\models T$, $p, q\in S_1(\emptyset)$ be non-algebraic
isolated types, $RC(p)=RC(q)=1$. Then $p(\mathcal{M})$ and
$q(\mathcal{M})$ are relatively $\mathcal{H}$-independent
$\Leftrightarrow$ $p\perp^w q$.
\end{lemma}

Proof of Lemma \ref{l_ind}. $(\Rightarrow)$ Let $p(\mathcal{M})$
and $q(\mathcal{M})$ are relatively $\mathcal{H}$-independent.
Assume the contrary: $p\not\perp^w q$. By quite o-minimality there
is an $\emptyset$-definable bijection $f: p(\mathcal{M})\to
q(\mathcal{M})$. Since $RC(p)=RC(q)=1$, this bijection is strictly
monotonic. Take an arbitrary open convex set $A'\subseteq
p(\mathcal{M})$ and consider $f(A')$. By strict monotonicity of
$f$ the image $f(A')$ is also open convex set. Take arbitrary $a,b
\in f(A')$ with the condition $a<b$. Then let $B':=\{c\in
q(\mathcal{M}) \mid a<c<b\}$. Then there is no $\mathcal{M}_1\prec
\mathcal{M}$ such that $A'=p(\mathcal{M})\cap M_1$ and
$B'=q(\mathcal{M})\cap M_1$.

$(\Leftarrow)$ Let $p\perp^w q$. Then by almost
$\omega$-categoricity of $T$ the family $\{p, q\}$ is orthogonal
over $\emptyset$, whence $p(\mathcal{M})$ and $q(\mathcal{M})$ are
relatively $\mathcal{H}$-independent. \endproof

\begin{corollary}
Let $\;T$ be an almost $\omega$-categorical quite o-minimal
theory, $\mathcal{M}\models T$, $p\in S_1(\emptyset)$ be a
non-algebraic isolated type, $RC(p)=n$, where $n>1$. Suppose that
$E_1(x,y)$, $E_2(x,y)$, $\ldots$, $E_{n-1}(x,y)$ are
$\emptyset$-definable equivalence relations partitioning
$p(\mathcal{M})$ into infinitely many infinite convex classes so
that $E_1(a,\mathcal{M})\subset
E_2(a,\mathcal{M})\subset\ldots\subset E_{n-1}(a,\mathcal{M})$ for
any $a\in p(\mathcal{M})$. Then the following holds:

$1)$ Every $E_1$-class is relatively $\mathcal{H}$-free;

$2)$ Any two $E_1$-classes are relatively $\mathcal{H}$-independent;

$3)$ For any $2\le i\le n-1$ every $E_i$-class is not relatively
$\mathcal{H}$-free.
\end{corollary}

\begin{example}\rm
Let $\mathcal{M}=\langle \mathbb{Q}\times\mathbb{Q};
<,E^2,f^1\rangle$ be a linearly ordered structure, where
$\mathbb{Q}\times\mathbb{Q}$ is lexicographically ordered. The
symbol $E$ is interpreted by a binary relation defined as follows:
$E(a,b) \Leftrightarrow n_1=n_2$ for any $a=(n_1, m_1), b=(n_2,
m_2)\in \mathbb{Q}\times\mathbb{Q}$. The symbol $f$ is interpreted
by a unary function defined by the equality $f((n,m))=(n+1, m)$
for all $(n,m)\in \mathbb{Q}\times\mathbb{Q}$.

Obviously,  $E(x,y)$ is an equivalence relation partitioning $M$
into infinitely many infinite convex classes.

It can be established that $Th(\mathcal{M})$ is a quite o-minimal
theory, and it is not almost $\omega$-categorical. Note that
$p(x):=\{x=x\}\in S_1(\emptyset)$ is a non-algebraic isolated
type, $RC(p)=2$, every $E$-class is relatively $\mathcal{H}$-free,
however $E(a,\mathcal{M})$ and $E(f(a), \mathcal{M})$ are not
relatively $\mathcal{H}$-independent for each $a\in M$.
\end{example}

\begin{definition}
{\rm Let $T$ be a weakly o-minimal theory, $\mathcal{M}\models T$,
$p\in S_1(\emptyset)$ be a non-algebraic isolated type. Let
$E(x,y)$ be an $\emptyset$-definable equivalence relation
partitioning $p(\mathcal{M})$ into in\-fi\-ni\-te\-ly many
infinite convex classes. If $A\subseteq p(\mathcal{M})$ then we
denote by $A/E$ the set of representatives of $E$-classes having a
nonempty intersection with $A$. We say that $p(\mathcal{M})$ is
{\it relatively $(\mathcal{H}, E)$-free} if for any convex
$A'\subseteq p(\mathcal{M})$ such that $A'/E$ is a open set the
equality $A'=p(\mathcal{M})\cap M_1$ holds for some $M_1\in
H(\mathcal{M})$.}
\end{definition}

Note that in the latter definition in case of dense ordering of
$p(\mathcal{M})$ the convexity of $A'$ is essential. Indeed, let
$p(x):=\{U(x)\}$, $a_1, a_2\in p(\mathcal{M})$ such that
$\mathcal{M}\models E(a_1, a_2) \land a_1<a_2$. Consider the
following formula:
$$\phi(x,a_1, a_2):=U(x)\land [x\le a_1\lor x\ge a_2].$$

Let $A'=\phi(\mathcal{M},a_1,a_2)$. Obviously, $A'\subseteq
p(\mathcal{M})$, $A'$ is not convex, $A'/E$ is open convex set,
but there is no $\mathcal{M}_1\prec \mathcal{M}$ such that $A'=
p(\mathcal{M})\cap M_1$.

\begin{proposition} \label{pr_rf}
Let $T$ be an almost $\omega$-categorical quite o-minimal theory,
$\mathcal{M}\models T$, $p\in S_1(\emptyset)$ be a non-algebraic
isolated type, $E(x,y)$ be an $\emptyset$-definable equivalence
relation par\-ti\-tio\-ning $p(\mathcal{M})$ into infinitely many
infinite convex classes. Then $p(\mathcal{M})$ is relatively
$(\mathcal{H},E)$-free $\Leftrightarrow$ for any
$\emptyset$-definable equivalence relation $E'(x,y)$ partitioning
$p(\mathcal{M})$ into infinitely many infinite convex classes we
have $E'(a,\mathcal{M})\subseteq E(a,\mathcal{M})$ for some $a\in
p(\mathcal{M})$.
\end{proposition}

Proof of Proposition \ref{pr_rf}. $(\Rightarrow)$ Let
$p(\mathcal{M})$ be relatively $(\mathcal{H}, E)$-free. By almost
$\omega$-ca\-te\-go\-ri\-ci\-ty of $T$ $RC(p)<\omega$, i.e. there
is an $\emptyset$-definable equivalence relation $E^*(x,y)$
partitioning $p(\mathcal{M})$ into infinitely many infinite convex
classes, and for any $\emptyset$-definable equivalence relation
$E'(x,y)$ partitioning $p(\mathcal{M})$ into infinitely many
convex classes we have $E'(a,\mathcal{M})\subseteq
E^*(a,\mathcal{M})$ for some $a\in p(\mathcal{M})$. Then we assert
that $E(x,y)\equiv E^*(x,y)$ on $p(\mathcal{M})$. If this is not
true then either $E(a,\mathcal{M})\subset E^*(a,\mathcal{M})$ or
$E^*(a,\mathcal{M})\subset E(a,\mathcal{M})$. In the first case we
have a contradiction with relative $(\mathcal{H}, E)$-freedom of
$p(\mathcal{M})$. In the second case we have a contradiction with
that $E^*(x,y)$ is the greatest.

$(\Leftarrow)$ Let $E(x,y)$ be the greatest $\emptyset$-definable
equivalence relation partitioning $p(\mathcal{M})$ into infinitely
many convex classes. Then if as $A'$ we take an arbitrary convex
subset of $p(\mathcal{M})$ so that  $A'/E$ is open then we easily
find $\mathcal{M}_1\prec \mathcal{M}$ with the condition
$A'=p(\mathcal{M})\cap M_1$.
\endproof

\begin{example}\rm
Let $\mathcal{M}=\langle M,<, E^2_i\rangle_{i\in \omega}$ be a
linearly ordered structure, where for each $i\in\omega$ $E_i(x,y)$
defines an equivalence relation partitioning $M$ into infinitely
many convex classes, and $E_i$ partitions every $E_{i+1}$-class
into infinitely many $E_i$-classes, every $E_i$-class is convex
and open so that $E_i$-subclasses of each  $E_{i+1}$-class are
densely ordered without endpoints.

It can be established that $Th(\mathcal{M})$ is a quite o-minimal
theory non-being almost $\omega$-ca\-te\-go\-ri\-cal. Obviously,
$M$ is 1-indiscernible, i.e. $p(x):=\{x=x\}\in S_1(\emptyset)$. It
is not difficult to see that  $p(\mathcal{M})$ is relatively
$(\mathcal{H}, E_i)$-free for any $i\in\omega$.
\end{example}

\begin{definition}
{\rm Let $T$ be a weakly o-minimal theory, $\mathcal{M}\models T$,
$p_1,p_2\in S_1(\emptyset)$ be non-algebraic isolated types. Let
$E_1(x,y), E_2(x,y)$ be $\emptyset$-definable equivalence
relations partitioning $p_1(\mathcal{M})$ and $p_2(\mathcal{M})$
respectively into infinitely many convex classes. Suppose that
$p_1(\mathcal{M})$ is relatively $(\mathcal{H}, E_1)$-free and
$p_2(\mathcal{M})$ is relatively $(\mathcal{H}, E_2)$-free. We say
that $p_1(\mathcal{M})$ and $p_2(\mathcal{M})$ are {\it relatively
$(\mathcal{H}, E_1, E_2)$-independent} if for any convex
$A'\subseteq p_1(\mathcal{M})$ and $B'\subseteq p_2(\mathcal{M})$
such that $A'/E_1$ and $B'/E_2$ are open sets there is $M_1\in
H(\mathcal{M})$ such that $A'=p_1(\mathcal{M})\cap M_1$ and
$B'=p_2(\mathcal{M})\cap M_1$.}
\end{definition}

\begin{proposition}\label{pr_ri}
Let $T$ be an almost $\omega$-categorical quite o-minimal theory,
$\mathcal{M}\models T$, $p_1, p_2\in S_1(\emptyset)$ be
non-algebraic isolated types. Let $E_1(x,y)$, $E_2(x,y)$ be
$\emptyset$-definable equivalence relations partitioning
$p_1(\mathcal{M})$ and $p_2(\mathcal{M})$ respectively into
infinitely many convex classes. Suppose that $p_1(\mathcal{M})$ is
relatively $(\mathcal{H}, E_1)$-free, $p_2(\mathcal{M})$ is
relatively $(\mathcal{H}, E_2)$-free. Then $p_1(\mathcal{M})$ and
$p_2(\mathcal{M})$ are {\it relatively $(\mathcal{H}, E_1,
E_2)$-independent} $\Leftrightarrow$ $p_1\perp^w p_2$.
\end{proposition}

Proof of Proposition \ref{pr_ri}. Let $p_1(\mathcal{M})$ and
$p_2(\mathcal{M})$ be relatively $(\mathcal{H}, E_1,
E_2)$-independent. Assume the contrary: $p_1\not\perp^w p_2$. By
quite o-minimality there is an $\emptyset$-definable bijection $f:
p_1(\mathcal{M})\to p_2(\mathcal{M})$, whence $RC(p_1)=RC(p_2)$
and $f(E_1(a,\mathcal{M}))=E_2(f(a),\mathcal{M})$ for any $a\in
p_1(\mathcal{M})$. Take an arbitrary convex set $A'\subseteq
p_1(\mathcal{M})$ with open $A'/E_1$ and consider $f(A')$.
Obviously, $f(A')$ is convex and $f(A')/E_2$ is open. Take
arbitrary $E_2$-classes $C=E_2(a,\mathcal{M})$ and
$D=E_2(b,\mathcal{M})$ for some $a,b\in p_2(\mathcal{M})$ with the
condition $C<D$ lying in $f(A')$. Then let $B':=\{e\in
p_2(\mathcal{M})\mid C<e<D\}$. Obviously, $B'$ will be also
convex, and $B'/E_2$ will be open. It is easily to see that there
is no $\mathcal{M}_1\prec \mathcal{M}$ such that
$A'=p_1(\mathcal{M})\cap M_1$ and $B'=p_2(\mathcal{M})\cap M_1$.
\endproof

\begin{corollary}
Let $T$ be an almost $\omega$-categorical quite o-minimal theory,
$\mathcal{M}\models T$, $p_1, p_2\in S_1(\emptyset)$ be
non-algebraic isolated types, and suppose that there exists an
$\emptyset$-definable bijection $f: p_1(\mathcal{M})\to
p_2(\mathcal{M})$. Let $E_1(x,y)$ be an $\emptyset$-definable
equivalence relation partitioning $p_1(\mathcal{M})$ into
infinitely many convex classes. Define on the set
$p_2(\mathcal{M})$ the relation $E_2(x,y)$ as follows:
$$ \mbox{for any } a, b \in p_2(\mathcal{M}) \quad E_2(a,b) \Leftrightarrow E_1(f^{-1}(a), f^{-1}(b)).$$
Then  $p_1(\mathcal{M})$ is relatively $(\mathcal{H}, E_1)$-free
$\Leftrightarrow$ $p_2(\mathcal{M})$ is relatively $(\mathcal{H},
E_2)$-free.
\end{corollary}

\vskip 5mm
Further we extend definitions of relative
$\mathcal{H}$-freedom, relative $\mathcal{H}$-independence,
relative $(\mathcal{H}, E)$-freedom and relative $(\mathcal{H},
E_1, E_2)$-independence on non-isolated 1-types.

Recall that if  $A$ is an arbitrary subset of a linearly ordered
structure $\mathcal{M}$ then we denote by $A^+$ (and respectively
by $A^-$) the sets of elements $b$ of the considered structure
with the condition $A<b$ ($b<A$).

\begin{definition} \rm\cite{bbs1}
Let $\mathcal{M}$ be a weakly o-minimal structure, $A\subseteq M$,
$p\in S_1(A)$ be a non-algebraic type. We say that $p$ is {\it
quasirational to right (left)} if there is an $A$-definable convex
formula $U_p(x)\in p$ such that for any sufficiently saturated
model $\mathcal{N}\succ \mathcal{M}$,
$U_p(\mathcal{N})^+=p(\mathcal{N})^+$
($U_p(\mathcal{N})^-=p(\mathcal{N})^-$). A non-isolated 1-type is
called {\it quasirational} if it is either quasirational to right
or quasirational to left. A non-quasirational non-isolated 1-type
is called {\it irrational}.
\end{definition}

Obviously, an 1-type being simultaneously quasirational to right
and quasirational to left is isolated.

We say that a convex set $A$ is {\it open to right (left)} if
there is $a\in A$ such that for any $b>a$ ($b<a$) there exists a
neighbourhood of the element  $b$ containing in $A$. Obviously, a
set being simultaneously open to right and open to left is open.

\begin{definition}\label{def_hfree}
\rm Let $T$ be a weakly o-minimal theory, $\mathcal{M}$ be a
sufficiently saturated model for $T$, $p\in S_1(\emptyset)$ be a
non-isolated type, $RC(p)=1$. If $p$ is quasirational to right
(left) then we say $p(\mathcal{M})$ is {\it relatively
$\mathcal{H}$-free} if for any open to right (left) convex
$A'\subseteq p(\mathcal{M})$ the equality $A'=p(\mathcal{M})\cap
M_1$ holds for some $M_1\in H(\mathcal{M})$. If $p$ is irrational
then it is sufficiently to take an arbitrary convex set as $A'$.
\end{definition}

\begin{lemma}\label{l_free_ni}
Let $T$ be an almost $\omega$-categorical quite o-minimal theory,
$\mathcal{M}$ be a sufficiently saturated model for $T$, $p\in
S_1(\emptyset)$ be a non-isolated type. Then $p(\mathcal{M})$ is
relatively $\mathcal{H}$-free $\Leftrightarrow$ $RC(p)=1$.
\end{lemma}

Proof of Lemma \ref{l_free_ni}. $(\Rightarrow)$ Indeed, if
$RC(p)>1$ then there is an $\emptyset$-definable equivalence
relation $E(x,y)$ partitioning $p(\mathcal{M})$ into infinitely
many infinite convex classes. Obviously, there is no
$\mathcal{M}_1\prec M$ such that
$E(a,\mathcal{M})=p(\mathcal{M})\cap M_1$ for some $a\in
p(\mathcal{M})$.

$(\Leftarrow)$ If $RC(p)=1$ then by analogy with proof of Lemma
\ref{l_free} it is established that $p(\mathcal{M})$ is
indiscernible over $\emptyset$. Then if $p$ is quasirational to
right (left) then for any open to right (left) convex set
$A'\subseteq p(\mathcal{M})$ there is $\mathcal{M}_1\prec
\mathcal{M}$ such that $A'=p(\mathcal{M})\cap M_1$. If $p$ is
irrational then for any convex set $A'\subseteq p(\mathcal{M})$
(including the case when $A'=\{a\}$ for some $a\in
p(\mathcal{M})$) there exists $\mathcal{M}_1\prec \mathcal{M}$
with $A'=p(\mathcal{M})\cap M_1$.
\endproof

\begin{definition}
{\rm Let $T$ be a weakly o-minimal theory, $\mathcal{M}$ be a
sufficiently saturated model for $T$, $p, q\in S_1(\emptyset)$ be
non-isolated types, $RC(p)=RC(q)=1$. We say that $p(\mathcal{M})$
and $q(\mathcal{M})$ {\it relatively $\mathcal{H}$-independent} if
for any convex sets $A'\subseteq p(\mathcal{M})$ and $B'\subseteq
q(\mathcal{M})$ corresponding to $p$ and $q$ (as in Definition
\ref{def_hfree}) there exists $M_1\in H(\mathcal{M})$ such that
$A'=p(\mathcal{M})\cap M_1$ and $B'=q(\mathcal{M})\cap M_1$.}
\end{definition}

\begin{proposition}\label{pr_norg} {\rm \cite{bbs1}}
Let $T$ be a weakly o-minimal theory, $\mathcal{M}\models T$,
$A\subseteq M$, $p, q\in S_1(A)$ be non-algebraic types,
$p\not\perp^w q$. Then:

$(1)$ $p$ is irrational $\Leftrightarrow$ $q$ is irrational;

$(2)$ $p$ is quasirational $\Leftrightarrow$ $q$ is quasirational.
\end{proposition}

\begin{lemma}\label{l_ind_ni}
Let $T$ be an almost $\omega$--categorical quite o-minimal theory,
$\mathcal{M}$ be a sufficiently saturated model for $T$, $p, q\in
S_1(\emptyset)$ be non-isolated types, $RC(p)=RC(q)=1$. Then
$p(\mathcal{M})$ and $q(\mathcal{M})$ are relatively
$\mathcal{H}$-independent $\Leftrightarrow$ $p\perp^w q$.
\end{lemma}

Proof of Lemma \ref{l_ind_ni}. If  $p\not\perp^w q$ then by
Proposition \ref{pr_norg} $p$ and $q$ are simultaneously either
quasirational or irrational. Without loss of generality, suppose
that $p$ and $q$ are quasirational. By quite o-minimality there is
an $\emptyset$-definable bijection $f: p(\mathcal{M})\to
q(\mathcal{M})$. Since the convexity ranks of the types are equal
to 1 then this bijection is strictly monotonic. For definiteness,
suppose that $p$ is quasirational to right. Then if $f$ is
strictly increasing (decreasing) then $q$ will be quasirational to
right (left). Take an arbitrary open to right convex set
$A'\subseteq p(\mathcal{M})$ and consider $f(A')$. If $f$ is
strictly increasing (decreasing) then $f(A')$ will be also open to
right (left) convex set. Take arbitrary $a,b \in f(A')$ with
$a<b$. Then let $B':=\{c\in q(\mathcal{M}) \mid a<c<b\}$. Then
there is no $\mathcal{M}_1\prec \mathcal{M}$ such that
$A'=p(\mathcal{M})\cap M_1$ and $B'=q(\mathcal{M})\cap M_1$.
\endproof

\begin{definition}\label{def_hind}
{\rm Let $T$ be a weakly o-minimal theory, $\mathcal{M}$ be a
sufficiently saturated model for $T$, $p\in S_1(\emptyset)$ be a
non-isolated type. Let $E(x,y)$ be an $\emptyset$-definable
equivalence relation partitioning $p(\mathcal{M})$ into infinitely
many infinite convex classes. If $p$ is quasirational to right
(left) then we say $p(\mathcal{M})$ is {\it relatively
$(\mathcal{H},E)$-free} if for any convex $A'\subseteq
p(\mathcal{M})$ such that  $A'/E$ is open to right (left) set the
equality $A'=p(\mathcal{M})\cap M_1$ holds for some $M_1\in
H(\mathcal{M})$. If $p$ is irrational then it is sufficiently to
take any open convex subset of $p(\mathcal{M})$ as $A'$, leaving
the type of the set $A'/E$ for an arbitrariness.  }
\end{definition}

\begin{proposition} \label{pr_rf_ni}
Let $T$ be an almost $\omega$-categorical quite o-minimal theory,
$\mathcal{M}$ be a suf\-fi\-cient\-ly saturated model for $T$,
$p\in S_1(\emptyset)$ be a non-isolated type, $E(x,y)$ be an
$\emptyset$-definable equivalence relation partitioning
$p(\mathcal{M})$ into infinitely many infinite convex classes.
Then $p(\mathcal{M})$ is relatively $(\mathcal{H},E)$-free
$\Leftrightarrow$ $E(x,y)$ is the greatest $\emptyset$-definable
equivalence relation partitioning $p(\mathcal{M})$ into infinitely
many convex classes.
\end{proposition}

Proof of Proposition \ref{pr_rf_ni}. $(\Rightarrow)$ is proved
similarly Proposition \ref{pr_rf}.

$(\Leftarrow)$  If $p$ is quasirational to right (left) then
taking as $A'$ an arbitrary convex subset of $p(\mathcal{M})$ with
the condition that $A'/E$ is open to right (left) we easily find
$\mathcal{M}_1\prec \mathcal{M}$ such that $A'=p(\mathcal{M})\cap
M_1$. If $p$ is irrational then take as $A'$ an arbitrary open
convex subset of $p(\mathcal{M})$.
\endproof

\begin{definition}
{\rm Let $T$ be a weakly o-minimal theory, $\mathcal{M}$ be a
sufficiently saturated model for $T$, $p_1,p_2\in S_1(\emptyset)$
be non-isolated types. Let $E_1(x,y), E_2(x,y)$ be
$\emptyset$-definable equivalence relations partitioning
$p_1(\mathcal{M})$ and $p_2(\mathcal{M})$ respectively into
infinitely many infinite convex classes. Suppose that
$p_1(\mathcal{M})$ is relatively $(\mathcal{H}, E_1)$-free and
$p_2(\mathcal{M})$ is relatively $(\mathcal{H}, E_2)$-free. We say
that $p_1(\mathcal{M})$ and $p_2(\mathcal{M})$ are {\it relatively
$(\mathcal{H}, E_1, E_2)$-independent} if for any convex
$A'\subseteq p_1(\mathcal{M})$ and $B'\subseteq p_2(\mathcal{M})$
corresponding to $p_1$ and $p_2$ (as in Definition \ref{def_hind})
there is $M_1\in H(\mathcal{M})$ such that
$A'=p_1(\mathcal{M})\cap M_1$ and $B'=p_2(\mathcal{M})\cap M_1$.}
\end{definition}

\begin{proposition}\label{pr_ri_ni}
Let $T$ be an almost $\omega$-categorical quite o-minimal theory,
$\mathcal{M}$ be a suf\-fi\-cient\-ly saturated model for $T$,
$p_1, p_2\in S_1(\emptyset)$ be non-isolated types. Suppose that
$p_1(\mathcal{M})$ is relatively $(\mathcal{H}, E_1)$-free and
$p_2(\mathcal{M})$ is relatively $(\mathcal{H}, E_2)$-free, where
$E_1(x,y)$, $E_2(x,y)$ are $\emptyset$-definable equivalence
relations partitioning $p_1(\mathcal{M})$ and $p_2(\mathcal{M})$
respectively into infinitely many infinite convex classes. Then
$p_1(\mathcal{M})$ and $p_2(\mathcal{M})$ are {\it relatively
$(\mathcal{H}, E_1, E_2)$-independent} $\Leftrightarrow$
$p_1\perp^w p_2$.
\end{proposition}

Proof of Proposition \ref{pr_ri_ni}. If $p_1\not\perp^w p_2$ then
by Proposition \ref{pr_norg} $p_1$ and $p_2$ are
si\-mul\-ta\-neous\-ly either quasirational or irrational. Without
loss of generality, suppose that $p_1$ and $p_2$ are
quasirational. By quite o-minimality there is an
$\emptyset$-definable bijection $f: p_1(\mathcal{M})\to
p_2(\mathcal{M})$. For definiteness, let $p_1$ be quasirational to
right. Then take an arbitrary convex $A'\subseteq
p_1(\mathcal{M})$ with the condition that $A'/E_1$ is open to
right. Obviously, $f(A')$ will be convex. If $f$ is strictly
increasing (decreasing) on $p_1(\mathcal{M})/E_1$ then $f(A')/E_2$
will be open to right (left). Taking arbitrary $E_2$-classes
$E_2(a,\mathcal{M})$ and $E_2(b,\mathcal{M})$ for some $a,b\in
p_2(\mathcal{M})$ with $E_2(a,\mathcal{M})<E_2(b,\mathcal{M})$
lying in $f(A')$, and considering $B':=\{h\in p_2(\mathcal{M})\mid
E_2(a,\mathcal{M})<h<E_2(b,\mathcal{M})\}$, we see that $B'$ is
convex, and $B'/E_2$ is open Obviously, there is no
$\mathcal{M}_1\prec \mathcal{M}$ such that
$A'=p_1(\mathcal{M})\cap M_1$ and $B'=p_2(\mathcal{M})\cap M_1$.
\endproof

\section{On freedom and independence in hypergraphs of models of theories of unars}
\noindent

In the case of unary theories, if sets $A_i=f^{-k_i}(a_i)$,
$k_i>0$, are $\mathcal{H}$-independent, then the restriction of
$\mathcal{H}$ on $\bigcup\limits_i A_i$ is represented in the form
of disjoint complete union of restrictions
$\mathcal{H}_i=\mathcal{H}\upharpoonright A_i$.

\begin{example} \rm Consider a connected free unar $\mathcal{M}=\langle
M,f\rangle$, i.e. a connected unar non-having cycles and such that
every element has infinitely many $f$-preimages. Consider also a
hypergraph $\mathcal{H}$ of elementary subsystems of
$\mathcal{M}$. Then for every element $a\in M$ and pairwise
distinct elements $a_i\in f^{-k}(a)$, $k>0$, the sets
$A_i=f^{-k_i}(a_i)$ are $\mathcal{H}$-independent for any $k_i>0$.
The restriction of $\mathcal{H}$ on $\bigcup\limits_i A_i$ is
represented in the form of disjoint complete union of restrictions
$\mathcal{H}_i=\mathcal{H}\upharpoonright A_i$.

On the other hand, the sets $f^{-k}(a)$ and $f^{-m}(b)$ for $k,m>0$
and $b\in f^{-k}(a)$ are not $\mathcal{H}$-independent, since $b\not\in Z\in
H(\mathcal{M})$ implies $f^{-m}(b)\cap Z=\emptyset$.

Thus, $\mathcal{H}$-independence of the sets $f^{-k}(a)$ and
$f^{-m}(b)$ is equivalent to their disjointness, and also the
conditions $b\notin\bigtriangleup_f(a)$ and
$a\notin\bigtriangleup_f(b)$, where
$\bigtriangleup_f(a)=\bigcup\limits_{n\in\omega}f^{-n}(a)$ is a
lower cone of the root $a$ \cite{CCMCT14}. The indicated
description of $\mathcal{H}$-in\-de\-pen\-den\-ce is naturally is
extended on an arbitrary family of the sets $f^{-k_i}(a_i)$.
\end{example}

By Properties 1 and 3, of $\mathcal{H}$-freedom and
$\mathcal{H}$-in\-de\-pen\-den\-ce, in general case the
$\mathcal{H}$-in\-de\-pen\-den\-ce of the sets $A_i=f^{-k_i}(a_i)$
for $a_i\in f^{-k}(a)$, allowing to conduct an indicated
decomposition, implies their infinity and presence of completeness
of the types $p_i(x)$ over $\{a_i\}$ isolated by the formulas
$f^{k_i}(x)\approx a_i$.

If the set $A_i=f^{-k_i}(a_i)$ is finite then its inclusion in
$Z\in H(\mathcal{M})$ is equivalent to inclusion of $a_i$ in $Z$.

\end{document}